\documentclass[a4paper,11pt,twoside,english]{article}
\usepackage{babel}
\usepackage{latexsym,amsfonts}
\usepackage{amsthm}
\usepackage{amsmath}
\usepackage{mathtools}
\usepackage{paralist}

\usepackage{lastpage}
\usepackage{fancyhdr}
\fancypagestyle{plain}{

  \fancyhead[L,C]{}
  \fancyhead[R]{\today}
  \fancyfoot[R,L]{}
  \fancyfoot[C]{\thepage \ of \pageref{LastPage}}
}

\title{Higher derivatives of the inverse tangent function and a summation formula involving binomial coefficients}
\author{Jan-David Hardtke}
\date{}

\providecommand{\sm}{\setminus}
\providecommand{\ssq}{\subseteq}
\providecommand{\N}{\ensuremath{\mathbb{N}}}
\providecommand{\R}{\ensuremath{\mathbb{R}}}

\providecommand{\keywords}[1]{
{\let\thefootnote=\relax
\footnote{{\em Keywords}: #1}}
\addtocounter{footnote}{-1}
}

\providecommand{\AMS}[1]{
{\let\thefootnote=\relax
\footnote{{\em AMS Subject Classification} (2010): #1}}
\addtocounter{footnote}{-1}
}

\providecommand{\address}{
{\sc \noindent Department of Mathematics \\
Universit\"at Leipzig\\
Augustusplatz 10, 04109 Leipzig \\
Germany \\}
}

\DeclarePairedDelimiter{\set}{\lbrace}{\rbrace}
\DeclarePairedDelimiter{\paren}{\lparen}{\rparen}

\DeclarePairedDelimiter{\floor}{\lfloor}{\rfloor}

\theoremstyle{definition}
\newtheorem{definition}{Definition}[section]
\newtheorem*{definition*}{Definition}
\theoremstyle{plain}

\newtheorem*{lemma*}{Lemma}
\newtheorem{proposition}[definition]{Proposition}
\newtheorem*{proposition*}{Proposition}

\newtheorem*{theorem*}{Theorem}
\newtheorem{corollary}[definition]{Corollary}
\newtheorem*{corolary*}{Corollary}

\newenvironment{Proof}[1][\proofname]{\begin{proof}[#1] \setlength{\parindent}{0pt}}{\end{proof}}
\newenvironment{Abstract}{\centering\begin{minipage}{0.8\textwidth} \noindent \small {\sc Abstract.}}{\end{minipage}\par}

\usepackage{color}
\definecolor{darkgreen}{rgb}{0,0.5,0}

\numberwithin{equation}{section}
\newtagform{colored}[\color{blue}]{\color{blue}(}{\color{blue})}
\usetagform{colored}

\usepackage[colorlinks,linkcolor=blue,citecolor=red,urlcolor=darkgreen]{hyperref}
\providecommand{\email}{{\it E-mail address:} \href{mailto:hardtke@math.uni-leipzig.de}{\tt hardtke@math.uni-leipzig.de}}

\usepackage{amsrefs}

\begin{document}

\maketitle

\begin{Abstract}
In \cite{deiser}, O. Deiser and C. Lasser obtained an explicit formula for the 
$n$-th derivative of the inverse tangent function. We calculate this derivative 
by a different method based on Fa\`a di Bruno's formula. Comparing the two results 
leads to the following identity for binomial coefficients:
\begin{equation*}
\sum_{i=m}^{\floor{n/2}}\frac{(-1)^i}{4^i}\binom{i}{m}\binom{n-i}{i}
=\frac{(-1)^m}{2^n}\binom{n+1}{2m+1},
\end{equation*}
where $n,m\in \N_0$ and $m\leq \floor{n/2}$. As was pointed out to the author by C. Krattenthaler, this formula is a special case of Gau\ss's formula for the hypergeometric function $_2F_1$.
\end{Abstract}
\AMS{26C05; 26A09}

\section{Higher derivatives of arctan}\label{sec:arctan}
The following explicit formula for the $n$-th derivative of the inverse tangent function
$\arctan$ was proved by O. Deiser and C. Lasser in \cite{deiser}:
\begin{equation}\label{eq:arctan1}
\arctan^{(n)}(x)=(n-1)!\frac{q_{n-1}(x)}{(1+x^2)^n} \ \ \ \forall x\in \R,\,\forall n\in \N,
\end{equation}
where
\begin{equation}\label{eq:arctan2}
q_n(x):=(-1)^n\sum_{k\text{\,even},\,0\leq k\leq n}\binom{n+1}{k+1}(-1)^{k/2}x^{n-k} \ \ \ \forall x\in \R,\,\forall n\in \N_0
\end{equation}
(where $\N_0:=\N\cup \set*{0}$).\par 
Other representations of $\arctan^{(n)}$ can be found in \cites{deiser,qi} and references therein. Here we want to obtain yet another explicit expression for $\arctan^{(n)}$ by using Fa\`a di Bruno's formula. Comparing our result to \eqref{eq:arctan1} then leads to a summation formula for binomial coefficients (which is a special case of Gau\ss's formula for the hypergeometric function $_2F_1$, see section 2).\par 
We start by recalling Fa\`a di Bruno's formula for the $n$-th derivative of the composition of two functions (see for instance the survey article \cite{johnson} and references therein): Given two intervalls $I,J\ssq \R$ and two $n$-times differentiable functions
$f:J \rightarrow \R$ and $g:I \rightarrow J$, we put $h:=f\circ g$. Then $h$ is $n$-times differentiable and for every $x\in I$ we have
\begin{equation}\label{eq:faadibruno}
h^{(n)}(x)=\sum_{(l_1,\dots,l_n)\in T_n}\frac{n!}{l_1!l_2!\dots l_n!}
f^{(l_1+l_2+\dots+l_n)}(g(x))\prod_{i=1}^n\paren*{\frac{g^{(i)}(x)}{i!}}^{l_i},
\end{equation}
where $T_n:=\set*{(l_1,\dots,l_n)\in \N_0^n:\sum_{i=1}^nil_i=n}$.\par 
There is also a slightly different version of Fa\`a di Bruno's formula based on Bell polynomials. After the first version of this preprint was published on arxiv.org, the author discovered the paper \cite{qi}. In this work one can also find an explicit formula for
the higher derivatives of $\arctan$, the proof of which is based on the Bell polynomial version
of Fa\`a di Bruno's formula. But the formula for $\arctan^{(n)}$ that we will obtain below (Proposition \ref{prop:arctan}) is different from the one in \cite{qi}.\par 
Now, as a special case of \eqref{eq:faadibruno} one gets the following result.
\begin{proposition}\label{prop:specialchainrule}
Let $n\in \N$, $a>0$ and let $f:(0,\infty) \rightarrow \R$ be an $n$-times differentiable function. We put $h(x):=f(a+x^2)$ for all $x\in \R$. Then $h$ is $n$-times differentiable and
\begin{equation}\label{eq:specialchainrule}
h^{(n)}(x)=\sum_{k=0}^{\floor{n/2}}\frac{n!}{k!(n-2k)!}(2x)^{n-2k}f^{(n-k)}(a+x^2) \ \ \ \forall x\in \R.	
\end{equation}	
\end{proposition}

\begin{Proof}
We put $g(x):=a+x^2$ for $x\in \R$. Then $g^{\prime}(x)=2x$, $g^{\prime\prime}(x)=2$ and 
$g^{(k)}(x)=0$ for all $k\geq 3$ and all $x\in \R$.\par
Let $S:=\set*{(l_1,\dots,l_n)\in T_n:l_i=0\ \text{for}\ i=3,\dots,n}$. Since $h=f\circ g$ it follows from Fa\`a di Bruno's formula that
\begin{equation*}
h^{(n)}(x)=\sum_{(l_1,\dots,l_n)\in S}\frac{n!}{l_1!l_2!}f^{(l_1+l_2)}(a+x^2)(2x)^{l_1} \ \ \ \forall x\in \R,
\end{equation*}
which can be rewritten as \eqref{eq:specialchainrule}.
\end{Proof}

Alternatively, one can also prove this statement directly by induction, without using 
Fa\`a di Bruno's formula (see the Appendix).\par
Proposition \ref{prop:specialchainrule} now allows us to obtain the following formula for the higher derivatives of the inverse tangent function.
\begin{proposition}\label{prop:arctan}
For every $n\in \N_0$ and every $x\in \R$ one has	
\begin{equation*}
\arctan^{(n+1)}(x)=\frac{n!2^n(-1)^n}{(1+x^2)^{n+1}}\sum_{m=0}^{\floor{n/2}}a_{mn}x^{n-2m},
\end{equation*}
where
\begin{equation*}
a_{mn}:=\sum_{k=m}^{\floor{n/2}}\frac{(-1)^k}{4^k}\binom{k}{m}\binom{n-k}{k} \ \ \forall m=0,\dots,\floor{n/2}.
\end{equation*}
\end{proposition}

\begin{Proof}
Let $h(x):=1/(1+x^2)=\arctan^{\prime}(x)$ for all $x\in \R$ and $f(y):=1/y$ for all $y\in \R\sm\set*{0}$. Then $f^{(k)}(y)=k!(-1)^k/y^{k+1}$ holds for every $k\in \N_0$.\par
Since $h(x)=f(1+x^2)$ it follows from Proposition \ref{prop:specialchainrule} that
\begin{align*}
&\arctan^{(n+1)}(x)=h^{(n)}(x)=\sum_{k=0}^{\floor{n/2}}\frac{n!(-1)^{n-k}(n-k)!}{k!(n-2k)!(1+x^2)^{n+1-k}}(2x)^{n-2k}\\
&=\frac{n!2^n(-1)^n}{(1+x^2)^{n+1}}\sum_{k=0}^{\floor{n/2}}(1+x^2)^k\frac{(-1)^k(n-k)!}{k!(n-2k)!}\frac{x^{n-2k}}{4^k}\\
&=\frac{n!2^n(-1)^n}{(1+x^2)^{n+1}}\sum_{k=0}^{\floor{n/2}}(1+x^2)^k(-1)^k\binom{n-k}{k}\frac{x^{n-2k}}{4^k}.
\end{align*}
Using the binomial theorem, we obtain
\begin{align*}
&\arctan^{(n+1)}(x)=\frac{n!2^n(-1)^n}{(1+x^2)^{n+1}}\sum_{k=0}^{\floor{n/2}}\sum_{i=0}^k(-1)^k\binom{k}{i}\binom{n-k}{k}\frac{x^{n-2k+2i}}{4^k}\\
&=\frac{n!2^n(-1)^n}{(1+x^2)^{n+1}}\sum_{(k,i)\in A}(-1)^k\binom{k}{i}\binom{n-k}{k}\frac{x^{n-2k+2i}}{4^k},
\end{align*}
where $A:=\set*{(k,i):k\in \set*{0,\dots,\floor{n/2}},i\in \set*{0,\dots,k}}$.\par 
Let $A_m:=\set*{(k,i)\in A:k-i=m}$ for all $m=0,\dots,\floor{n/2}$. It follows that
\begin{equation*}
\arctan^{(n+1)}(x)=\frac{n!2^n(-1)^n}{(1+x^2)^{n+1}}\sum_{m=0}^{\floor{n/2}}x^{n-2m}\sum_{(k,i)\in A_m}\frac{(-1)^k}{4^k}\binom{k}{i}\binom{n-k}{k}.
\end{equation*}
But 
\begin{align*}
&\sum_{(k,i)\in A_m}\frac{(-1)^k}{4^k}\binom{k}{i}\binom{n-k}{k}=\sum_{k=m}^{\floor{n/2}}\frac{(-1)^k}{4^k}\binom{k}{k-m}\binom{n-k}{k}\\
&=\sum_{k=m}^{\floor{n/2}}\frac{(-1)^k}{4^k}\binom{k}{m}\binom{n-k}{k}=a_{mn} \ \ \ \forall m=0,\dots,\floor{n/2},
\end{align*}
which completes the proof.
\end{Proof}

\section{A summation formula involving binomial coefficients}
Now we can compare our result Proposition \ref{prop:arctan} with the simpler formula \eqref{eq:arctan1} that was found in \cite{deiser} and obtain the following summation formula
for binomial coefficients.
\begin{proposition}\label{prop:sums}
Let $n,m\in \N_0$ such that $m\leq \floor{n/2}$. Then we have
\begin{equation}\label{eq:binom}
\sum_{i=m}^{\floor{n/2}}\frac{(-1)^i}{4^i}\binom{i}{m}\binom{n-i}{i}
=\frac{(-1)^m}{2^n}\binom{n+1}{2m+1}.
\end{equation}
\end{proposition}

\begin{Proof}
From Proposition \ref{prop:arctan} and the result of \cite{deiser} (formulas \eqref{eq:arctan1} and \eqref{eq:arctan2}) it follows that
\begin{equation*}
2^n\sum_{m=0}^{\floor{n/2}}a_{mn}x^{n-2m}=\sum_{m=0}^{\floor{n/2}}\binom{n+1}{2m+1}(-1)^mx^{n-2m}
\end{equation*}
holds for all $x\in \R$. Hence
\begin{equation*}
\frac{(-1)^m}{2^n}\binom{n+1}{2m+1}=a_{mn}=\sum_{i=m}^{\floor{n/2}}\frac{(-1)^i}{4^i}\binom{i}{m}\binom{n-i}{i}
\end{equation*}
for every $m\in \set*{0,\dots,\floor{n/2}}$.
\end{Proof}

After the first version of this preprint was published on arxiv.org, the author was informed by Christian Krattenthaler that this formula is a special case of Gau\ss's formula for the hypergeometric function $_2F_1$. Here is a sketch of the argument: $_2F_1$ is defined by 
\begin{equation*}
_2F_1(a,b;c;z):=\sum_{k=0}^{\infty}\frac{(a)_k(b)_k}{(c)_k}\frac{z^k}{k!},
\end{equation*}
where $(q)_k$ is the Pochhammer symbol, i.\,e. $(q)_k:=q(q+1)\dots(q+k-1)$ for $k\geq 1$ and $(q)_0:=1$.\par 
Gau\ss's formula reads
\begin{equation*}
_2F_1(a,b;c;1)=\frac{\Gamma(c)\Gamma(c-a-b)}{\Gamma(c-a)\Gamma(c-b)},
\end{equation*}
where $\Gamma$ denotes the Gamma function (information on $\Gamma$ and $_2F_1$ can be found, for instance, in Chapter 1 resp. 2 of \cite{erdelyi}).\par
The sum in \eqref{eq:binom} can be expressed via $_2F_1$ as
\begin{align*}
&\sum_{i=m}^{\floor{n/2}}\frac{(-1)^i}{4^i}\binom{i}{m}\binom{n-i}{i}\\
&={_2F_1}(m-n/2,m-n/2+1/2;m-n;1)\frac{(-1)^m}{m!4^m} (n-2m+1)_m.
\end{align*}
Applying Gau\ss's formula then leads to the result of Proposition \ref{prop:sums}.\par 
\ \par
A consequence of \eqref{eq:binom} is the following formula (which is probably also known,
but the author was unable to find a reference).
\begin{corollary}\label{cor:sums}
For all $n\in \N_0$ we have
\begin{equation*}
\sum_{i=0}^n\frac{(-1)^i}{4^i(n+1-i)}\binom{2n+1-i}{i}=\begin{cases}
0\ \ \text{if}\ n\ \text{is\ odd},\\
\frac{4^{-n}}{n+1}\ \ \text{if}\ n\ \text{is\ even}.
\end{cases}
\end{equation*}
\end{corollary}

\begin{Proof}
For each $n\in \N_0$ we put
\begin{equation*}
C_n:=\sum_{i=0}^n\frac{(-1)^i}{4^i}\binom{2n-i}{i}.
\end{equation*}	
By Proposition \ref{prop:sums} (with $m=0$ and $2n$ instead of $n$) we have $C_n=4^{-n}(2n+1)$.
It follows that
\begin{equation}\label{eq:cor1}
C_{n+1}-\frac{C_n}{4}=\frac{2}{4^{n+1}}.
\end{equation}
On the other hand, we have
\begin{align*}
&C_{n+1}-\frac{C_n}{4}=\sum_{i=0}^{n+1}\frac{(-1)^i}{4^i}\binom{2n+2-i}{i}-\sum_{i=0}^n\frac{(-1)^i}{4^{i+1}}\binom{2n-i}{i}\\
&=1+\sum_{i=1}^{n+1}\frac{(-1)^i}{4^i}\paren*{\binom{2n+2-i}{i}+\binom{2n+1-i}{i-1}}.
\end{align*}
For $i\in \set*{1,\dots,n}$ we have
\begin{align*}
&\binom{2n+2-i}{i}+\binom{2n+1-i}{i-1}=\frac{(2n+2-i)!}{i!(2n+2-2i)!}+\frac{(2n+1-i)!}{(i-1)!(2n+2-2i)!}\\
&=\frac{(2n+2-i)!+i(2n+1-i)!}{i!(2n+2-2i)!}=\frac{(2n+1-i)!(2n+2)}{i!(2n+2-2i)!}\\
&=\frac{2n+2}{2n+2-2i}\binom{2n+1-i}{i}.
\end{align*}
It follows that
\begin{align*}
&\frac{1}{2n+2}\paren*{C_{n+1}-\frac{C_n}{4}}\\
&=\frac{1}{2n+2}\paren*{1+\sum_{i=1}^n\frac{(-1)^i(2n+2)}{4^i(2n+2-2i)}\binom{2n+1-i}{i}+\frac{(-1)^{n+1}2}{4^{n+1}}}\\
&=\frac{1}{2n+2}\paren*{1+\frac{(-1)^{n+1}2}{4^{n+1}}}+\sum_{i=1}^n\frac{(-1)^i}{4^i(2n+2-2i)}\binom{2n+1-i}{i}\\
&=\frac{(-1)^{n+1}}{4^{n+1}(n+1)}+\sum_{i=0}^n\frac{(-1)^i}{4^i(2n+2-2i)}\binom{2n+1-i}{i}.
\end{align*}
Together with \eqref{eq:cor1} this implies
\begin{equation*}
\frac{(-1)^{n+1}}{4^{n+1}(n+1)}+\sum_{i=0}^n\frac{(-1)^i}{4^i(2n+2-2i)}\binom{2n+1-i}{i}=
\frac{1}{4^{n+1}(n+1)}
\end{equation*}
and thus
\begin{equation*}
\sum_{i=0}^n\frac{(-1)^i}{4^i(n+1-i)}\binom{2n+1-i}{i}=\frac{2(1-(-1)^{n+1})}{4^{n+1}(n+1)}
=\begin{cases}
0\ \ \text{if}\ n\ \text{is\ odd},\\
\frac{4^{-n}}{n+1}\ \ \text{if}\ n\ \text{is\ even}.
\end{cases}
\end{equation*}
\end{Proof}

\section{Appendix}\label{sec:appendix}
Here we want to give a direct proof of Proposition \ref{prop:specialchainrule} via induction
(without using Fa\`a di Bruno's formula). We recall the Proposition's statement:\\
{\it 
If $n\in \N$, $a>0$, $f:(0,\infty) \rightarrow \R$ is an $n$-times differentiable function and $h(x):=f(a+x^2)$ for all $x\in \R$, then $h$ is $n$-times differentiable and
\begin{equation}\label{eq:appendix}
h^{(n)}(x)=\sum_{k=0}^{\floor{n/2}}\frac{n!}{k!(n-2k)!}(2x)^{n-2k}f^{(n-k)}(a+x^2) \ \ \ \forall x\in \R.	
\end{equation}	
		}

\begin{Proof}
For $n=1$ this follows immediately from the chain-rule.\par 
Now suppose that the statement is true for some $n\in \N$ and assume that $f$ is even $(n+1)$-times differentiable. Differentiating \eqref{eq:appendix} gives 
\begin{equation*}
h^{(n+1)}(x)=A(x)+B(x),
\end{equation*}
where
\begin{align*}
&A(x):=\sum_{k=0}^{\floor{n/2}}\frac{2^{n-2k}n!}{k!(n-2k)!}(n-2k)x^{n-2k-1}f^{(n-k)}(a+x^2),\\
&B(x):=\sum_{k=0}^{\floor{n/2}}\frac{n!}{k!(n-2k)!}(2x)^{n+1-2k}f^{(n+1-k)}(a+x^2).
\end{align*}
If $n$ is even, say $n=2l$, we get
\begin{equation*}
A(x)=\sum_{k=1}^l\frac{2^{n-2k+2}n!}{(k-1)!(n-2k+1)!}x^{n+1-2k}f^{(n+1-k)}(a+x^2)
\end{equation*}
and hence 
\begin{align*}
&h^{(n+1)}(x)-(2x)^{n+1}f^{(n+1)}(a+x^2)\\
&=\sum_{k=1}^l(2x)^{n+1-2k}f^{(n+1-k)}(a+x^2)n!\paren*{\frac{2}{(k-1)!(n-2k+1)!}+
\frac{1}{k!(n-2k)!}}\\
&=\sum_{k=1}^l\frac{(n+1)!}{k!(n+1-2k)!}(2x)^{n+1-2k}f^{(n+1-k)}(a+x^2).
\end{align*}
This implies
\begin{equation*}
h^{(n+1)}(x)=
\sum_{k=0}^{\floor{(n+1)/2}}\frac{(n+1)!}{k!(n+1-2k)!}(2x)^{n+1-2k}f^{(n+1-k)}(a+x^2).
\end{equation*}
If $n$ is odd ($n=2l+1$), then we have 
\begin{equation*}
A(x)=\sum_{k=1}^l\frac{2^{n-2k+2}n!}{(k-1)!(n-2k+1)!}x^{n+1-2k}f^{(n+1-k)}(a+x^2)+
\frac{2n!}{l!}f^{(l+1)}(a+x^2)
\end{equation*}
and by a similar calculation as in the even case we obtain
\begin{align*}
&h^{(n+1)}(x)-(2x)^{n+1}f^{(n+1)}(a+x^2)\\
&=\sum_{k=1}^l\frac{(n+1)!}{k!(n+1-2k)!}(2x)^{n+1-2k}f^{(n+1-k)}(a+x^2)+
\frac{2n!}{l!}f^{(l+1)}(a+x^2)\\
&=\sum_{k=1}^{l+1}\frac{(n+1)!}{k!(n+1-2k)!}(2x)^{n+1-2k}f^{(n+1-k)}(a+x^2).
\end{align*}
Thus we have again
\begin{equation*}
h^{(n+1)}(x)=
\sum_{k=0}^{\floor{(n+1)/2}}\frac{(n+1)!}{k!(n+1-2k)!}(2x)^{n+1-2k}f^{(n+1-k)}(a+x^2).
\end{equation*}
\end{Proof}
\ \\
{\bf Acknowledgements:} The author is grateful to Christian Krattenthaler for pointing out
to him the proof of \eqref{eq:binom} via Gau\ss's hypergeometric function.

\address
\email

\end{document}